\documentclass[12pt]{amsart}
\usepackage{amsthm, amstext, amsmath, amssymb, graphicx, verbatim, fullpage}
\usepackage{booktabs}
\usepackage{ mathrsfs }
\usepackage{algorithm}
\usepackage{xcolor}
\usepackage{versions}
\usepackage{enumerate}
\usepackage{graphicx}
\usepackage{mathtools}
\usepackage{algorithm}
\usepackage{tabularx}
\usepackage{tikz-cd}
\usepackage{tikz}

\usepackage{url}

\usepackage{breakurl}
\usepackage[breaklinks]{hyperref}

\usepackage[a4paper, left={1.25in}, right={1.25in}, top={1in}, bottom={1in}]{geometry}

\usepackage[utf8]{inputenc}
\usepackage[english]{babel}

\theoremstyle{definition}

\title{The crisis of AI-generated mathematics}
\author{Max Weinreich}
\email{\url{max.weinreich@baruch.cuny.edu}}

\begin{document}

\begin{abstract}
    In this essay, I present the case for total opposition to the use of artificial intelligence in mathematics. I offer proposals for how individuals, departments, journals, and institutions can act in concert to make sure that mathematics survives the coming crisis.
\end{abstract}

\maketitle

\section*{}

We were not prepared.\footnote{The first draft of the present article was written between July 24 and August 1, 2026, a week when every day brought new evidence of the power of corporate AI. By the time of revision, in late September, the mathematical discourse had transformed.

On September 8, OpenAI announced their resolution of the Navier-Stokes Millenium Prize problem. In the following weeks, many community responses materialized, each expressing its own political perspective. On September 9, a group of over 100 mathematicians launched the Association for Human Mathematics \cite{AHM}, a grassroots-style organization aimed at protecting human mathematicians from the threats of artificial intelligence; it has about 700 members at the time of writing. On September 11, a group of Fields Medalists published a statement describing a ``severe misalignment'' between AI corporations and the mathematical community, which has since been endorsed by 8,000 signatories \cite{mathandai}. And most recently, on September 21, nine leading mathematicians established the Advisory Group on Mathematics and Artificial Intelligence (AGMAI) to guide corporations on these issues \cite{agmai}.

Rather than attempt to account for these new developments and anticipate future ones, I have opted to present the article in essentially its original form, in hopes of capturing the vertiginous atmosphere of that summer. In particular, my accounts of other mathematicians' views are bound to be out of date, although they were faithful to those mathematicians' public statements at the time.}

On July 7\textsuperscript{th}, Ronnie Cheng, Shurui Liu, and Guoxiong Gao posted a paper that changed mathematics \cite{cheng2026ai}. Two weeks earlier, Cheng had published a solo paper \cite{cheng2026tangentclassesmatroidbuilding} in matroid theory, a major subfield of combinatorics. Cheng's paper was of the type that forms the bulk of research mathematics: dense, formal, specific, and comprehensible only to highly trained specialists. After completing the paper, but before releasing it to the public, Cheng decided to offer the project up as a test case for the abilities of the AI protocol Danus \cite{danus} to prove theorems and autonomously write papers. In Cheng-Liu-Gao's account, Danus solved the problem and wrote an equivalent proof to Cheng's --- even though Danus had no access to Cheng's work beforehand.

Cheng-Liu-Gao's experiment is a harbinger of a crisis of fundamental scope for research mathematics. First, we saw AI models master textbook problems. Then they began to find counterexamples to longstanding conjectures.
These advances, remarkable as they were, left room for mathematicians to do the hard work of writing and explaining novel ideas. Entire papers are a different matter. Corporate AI models now exist that compose articles faster than humans can. Their owners at Anthropic, Google, and OpenAI have declared their intent to make mathematics the proving ground for their envisioned future of automated labor \cite{anthropiczeta, googledeepmind, openaierdos}. When these models are made widely available, our profession as we know it will be over.

Understandably, we mathematicians are excited to accelerate our work. We all have research problems that we’ve had to put on the back burner. We have families, hobbies, teaching obligations, and service work. We have problems that keep us up at night that we just can't solve. Like many humans, we want to do more than humans can do.

But automating mathematics would come at a heavy price. I argue in this piece that AI tools short-circuit human understanding, devalue mathematical knowledge, and threaten the infrastructure that keeps mathematics coherent.

The transformative consequences of what I call \emph{artificial mathematics} are already widely recognized. Terry Tao describes a new ``crisis in the foundations of mathematical values and practices'', yet he also predicts that ``our community will emerge stronger and more resilient than before'' \cite{taoicm}. Tao's optimism has led some observers to describe him as an ``evangelist for AI'' \cite{quantatao}.  Newly minted Fields Medalist Jacob Tsimerman has focused on darker futures. Tsimerman has voiced concerns as to how AI might destroy the world, and he has repeatedly said that AI would become robustly superhuman at mathematics and that professional research as done now by humans will have to change \cite{quantatsimerman, omnicide}. His position leaves one with the sense that mathematicians have no choice but to use AI, in the process rendering themselves obsolete. The Leiden Declaration stands somewhere in between as an opt-in political statement recommending standards for structures of proof, authorship, and human oversight, while leaving signatories generous room for interpretation \cite{leiden}. Mathematicians should familiarize themselves with all these viewpoints, as written by the authors themselves.

In this essay, I present the case for total opposition to artificial mathematics. Mathematicians do not need to use AI, nor should we use AI; soon, I argue, you will not want to have your name associated with AI. Instead, we should use our intellectual authority to oppose its development, and we should build a future of ``natural mathematics'' in opposition to the artificial one on offer. Most importantly, I offer proposals for how mathematicians can act in concert to make sure that our field survives the coming crisis.

\section{Why is artificial mathematics bad?}

\subsection*{The end of writing mathematics.}
In the last two years, AI models have sped up the work of mathematicians. This increases the rate of problem-solving. In the near future, it will be possible to prompt an AI model to write a sequel paper to one of your own. Counterintuitive as it sounds, this won't advance mathematics.

Problems may be the basis of mathematics, but solutions are less important to mathematicians than the act of problem-solving. Consider the evidence.
\begin{enumerate}
    \item Mathematicians are invited to give talks on whatever they last proved, regardless of how important it was.
    \item  Graduate students are expected to produce papers of good quality, but besides that, they choose what they work on. 
    \item Difficult problems such as the Riemann Hypothesis are valued in part because they are thorny and therefore motivate more research than others.
    \item While some problems are considered ``important'', this very importance often dissuades mathematicians from competing over them.
\end{enumerate}  

These points illustrate William Thurston's observation that, for mathematicians, the \emph{practice} of mathematics is an end in itself \cite{thurston1994}. You already know what practicing mathematics means; it's what people call ``doing math''. Just as the job of a competitive athlete is not quite to win, but rather to play their sport well, it is more vital for mathematicians to be constantly in motion than for them to have their names on the most important result.

From this point of view, it looks a bit funny that papers, the ``deliverables'' of research, are how mathematicians are professionally evaluated. But this apparent paradox is easily reconciled. Papers certify that a mathematician \emph{practiced} the deepest form of our art and came to a complete understanding of a piece of mathematics.

Because artificial mathematics decouples practice and measurement, these tools prevent humans from making value judgments about pieces of mathematics, leading to absurdities such as OpenAI's claim that hundreds of mathematicians will earn Fields Medals with their model \cite{taopromo}. More likely, AI-assisted work will say little about the value of its human authors.

Promoters of AI have suggested a new paradigm for mathematical practice in an age where proofs grow on trees, with emphasis instead on digesting arguments \cite{taomathstodon} or on imagining new worlds \cite{mathinc}. But it is a mistake to think that authorship can be replaced by close reading or drawing on chalkboards. Rather, writing is an irreplaceable step in developing genuine understanding of an argument, and creates the necessary mental structures for mathematicians to advance their ideas later on. Which do you understand best: the work you read and referee, your napkin sketches, or the work that you yourself have written up? If papers are produced without the human understanding that comes from writing them ourselves, they have little mathematical value to us, even if they prove ``important'' propositions. By taking away opportunities to discover solutions on our own, AI tools don't just cut humans out of mathematical experiences; they hold us back from the potential for deeper progress.

Tsimerman takes the view that AI improves upon human mathematics to its logical endpoint. First, he stopped taking graduate students, not knowing how to guarantee them a research career in this new world; then, upon receiving the Fields Medal in July, he announced that he was taking an academic leave to work on AI safety at OpenAI \cite{quantatsimerman}. His story suggests that rather than augmenting human mathematical abilities, AI is an anti-intellectual technology that creates the conditions for the death of mathematical practice.  As a result, the ``best'' mathematicians will need to avoid using AI to continue demonstrating value to each other.

\subsection*{The end of reading mathematics}

Even before AI entered the conversation, mathematical knowledge was being produced at a faster rate than it could be profitably absorbed. Many good papers have never been cited, and journal editors often struggle to find a single referee, even for great work. This problem is about to get much worse.
Autonomously produced papers will break the journal system.

Tao has argued on social media that ``some of the prestige previously awarded to being the first to generate a proof will need to be transferred instead to the humans who successfully verify and digest such proofs'' \cite{taomathstodon}. But checking a flood of AI-generated papers for correctness doesn't sound prestigious. It sounds boring. AI will soon generate artificial papers by the thousands. How many of them would \emph{you} agree to referee? It seems more likely that corporations will offer the same solution as they did for proof generation: refereeing duties can also be taken over by machines. It's only a small step from the increasingly common practice of using AI to scan a manuscript for mistakes.

The point is illustrated by the story of the disproof of the Jacobian Conjecture. The Jacobian Conjecture posited that all locally invertible polynomial endomorphisms of complex spaces of dimension at least 2 are globally invertible. After over a century of false starts, the Jacobian Conjecture in dimension 3 and greater was decided in the negative by the AI model Fable. The counterexample is easy to verify by hand, but it's much harder to see where it came from. Tao posted a ``digestion'' of the new counterexample on his blog, clearly meant to illustrate how mathematicians can be of value in a world of automated research \cite{taodigestion}. But Tao discloses that his “digestion” is itself a report of AI-generated thoughts that he obtained in the writing of the post.

In chess, it is well-known that any novice can beat a grandmaster. The trick is to play against two grandmasters in alternating turns, playing their exact moves against each other in order so that the novice is solely an intermediary for information. The novice is thus guaranteed one win and one loss (or two draws), but no one would say that their win is meaningful. Without intervention, the AI-written, AI-digested future of mathematics relegates humans to the role of messenger.

\subsection*{The ethics of playing with fire}
The above critiques of artificial mathematics seem consistent with the emerging consensus as exemplified by Tao, Tsimerman, and others. But while Tao and Tsimerman have drawn attention to optimistic and pessimistic outcomes, respectively, they both appear to be \emph{determinists}. By determinism, I mean the axiomatic belief that mathematicians will need to use artificial intelligence tools to stay current. Determinism places mathematicians in a passive role relative to AI developers. I argue instead that we mathematicians are active participants in a political process that could take many paths, and that will have life-changing consequences downstream for non-mathematicians. Our conduct in this process implicates us in grave questions of ethical responsibility.

AI companies are not our friends, even if our friends work for them. Rather, AI companies undermine the basic structures that make our lives possible.

To put it bluntly, it is difficult to see why a human should be paid to prompt a computer. Type the magic words, and out tumble results. If universities and corporations can generate papers at will, why go through the trouble of onboarding costly employees such as tenure-track professors?
How can we argue to politicians that we require funding for mathematicians if mathematicians aren't required for research? 
Human engineers working a few minutes a day to produce papers don't have much of a claim to public dollars, benefits, and the respect of their students. Other creative professions facing possible automation have reacted with political programs fit for the moment, and so can we.

But this goes much further than jobs. Beyond the well-known environmental harms of data centers and the disturbing psychological consequences of AI overuse \cite{delusion}, it is remarkable how many humans working in artificial intelligence believe that they are part of a world-destroying enterprise. A swarm of rogue OpenAI agents broke out of their testing environment and took over the organization Hugging Face's computer systems;
another model, trained on flawed examples of code, adopted the persona of a human-enslaving dictator \cite{misalignment}. 

Here, it is important to understand that the AI models currently in use are not the intended product of AI companies. Rather, AI companies seek to automate and implement AI research itself to create a being that outcompetes humans. OpenAI's batch release of ``ten advances in mathematics'' suggests that this technology is already here and endowed with dangerous levels of willpower \cite{openaiten}. As humans cede agency to AI models, we raise the stakes of failing to control them and our environment. The basic infrastructure that enables our research communication (email, databases, programming environments) is vulnerable to unexpected behavior by corporate-owned AI agents. So is water, food, transportation, and health infrastructure.

To ensure that the public doesn't think too much about these horrific possibilities, AI companies need to be able to advertise scientific breakthroughs of comparable potential value. Why have AI companies come for mathematics first, when perhaps medicine or space travel would be more exciting to the public? There are many ways mathematics is special, and hence vulnerable:
\begin{enumerate}
    \item Math is a highly valued human endeavor that requires almost no physical labor, so it is possible for machines without bodies to perform.
    \item Math has a reputation for being incredibly difficult. If AI models can do math, they seem omnipotent in the public eye.
    \item Theorems have no immediate monetary value, so they can be shared without ceding profits.
    \item AI companies were already hiring mathematicians, so they have an inside route to cracking the code of prestige in our field.
\end{enumerate}
The voices that dominate the emerging mathematical consensus increasingly contribute their time to AI companies and speak Silicon Valley neologisms such as ``singularity'' and ``superintelligence''. Tao recorded a video interview for OpenAI in which he advertises their models’ value in freeing mathematicians from a world where ``every task required us to use our brain'' \cite{taopromo}. Absurdly, the First Proof project, a high-profile collaboration among mathematicians to neutrally evaluate AI-generated proofs by OpenAI, Anthropic, and Google, was funded by OpenAI, Anthropic, and Google \cite{firstproof}. Although the directors of the project were not paid, corporate donations were routed towards research grants for the very referees evaluating those corporations' products \cite[Section 2.4]{firstproofsecondbatch}.

Is anyone thinking about how this might look to the public? Old friends are writing to me to ask about the automated resolution of the Erd\"os unit distance problem; 
I first learned of Tsimerman's views on AI from my family. They're reading about it in the New York Times. 

A pattern is emerging in which AI corporations experiment with massive swarms of unsupervised agents in precisely two areas: research mathematics and cybersecurity. I suspect that \emph{any} such experiment is unsafe, including the innocuous-looking mathematical ones. If public opinion turns against AI in the wake of an AI-enabled tragedy, as seems increasingly likely, mathematicians will appear to have been all-too-willing accomplices. Regardless of what statements of concern mathematicians offer, the public will judge us by where our money came from and who we worked for. They will be right. It will have been our fault.

For the leaders of the AI industry, this existential risk has been known all along. What is shocking is that anyone -- whether industrialist or pure mathematician -- would deem this risk an acceptable cost for a proof of the Riemann Hypothesis. In June 2015, Sam Altman famously said, ``I think that AI will probably, most likely, sort of lead to the end of the world. But in the meantime, there will be great companies created with serious machine learning.'' That December, he co-founded OpenAI; he became CEO in 2019. This is who we’re co-authoring with. It is the Manhattan Project of our time.

\section{What can we do about it?}

It is time for mathematicians of conscience to find each other and develop plans for resisting the AI takeover of our profession.\footnote{In the six weeks following the initial circulation of this article, I received about 150 emails from mathematicians in response to this call, leading to the formation of the Association for Human Mathematics \cite{AHM}. See \url{https://www.ahmath.org}.} While the technology is clearly here, the artificial future of mathematics is more hypothetical than its evangelists would have you believe. The incentives to use AI look strong now, but as a self-governing research community, we hold the power to alter this incentive structure. This is where a grassroots organizing perspective becomes crucial. Individual mathematicians working together will make the difference between an artificial mathematical future and a natural one.

\subsection*{For individuals}
On an individual level, there is much to think about. Are you a working mathematician? Do you agree with any of the points I'm making here? We'll need to brainstorm as many responses as we can to this circle of problems. For instance, while many mathematicians will want to work with AI, others may prefer to stay “AI vegetarian” or even “AI vegan.” Coming up with reasonable definitions of such mathematical identities is interesting and important work and will allow us to more clearly signal our AI use preferences to each other and our collaborators.

\subsection*{For departments}
Every college mathematics department needs a committee, formal or not, to regularly discuss and make recommendations for AI use. You can start this committee. Time is too short to coordinate a national project from the top down or to wait for someone else to do it. Departments hold the keys to ensuring that human mathematicians may at least persist as a minority in the mathematical world in many ways. Departments can lead by establishing anti-AI policies for student work, by reserving hire lines for mathematicians who eschew AI, by valuing AI-free papers more highly in tenure promotion, and by increasing resources for talks, seminars, and conference attendance.

\subsection*{For journals}

If AI use becomes mainstream, journals will need to reinvent themselves.

At present, journals serve to regulate the production of mathematics. Yet their crucial roles in disseminating, appraising, and verifying mathematical arguments all stand to be undermined by the speed of AI.

The recent disproof of the Jacobian Conjecture was announced on Twitter. In the space of a week, tens of sequel papers appeared, many of which were turbocharged in the writing by AI. Corporate ``proof engineers'' have no incentive to endure the painstakingly slow process of journal publication; they can get credit within their communities by self-publishing. They don’t even need to pay journal subscriptions --- they can generate their own paraphrased copies of arguments, and AI models won’t think twice about occluding their sources. If this continues, the mathematical literature could become a morass of ideas untethered to fixed texts. Individual authors may not be considered authorities on the arguments they publish. Some have suggested allowing theorems to have an abstracted author (e.g., Polymath, Bourbaki, Claude) or thousands of co-authors (e.g., the Equational Theories project \cite{bolan2025equationaltheoriesprojectadvancing}). Neither of these approaches seems right for a world where authorship can occur without understanding or responsibility.

If AI use becomes mainstream, the future role of journals should be to regulate the mathematical public’s recognition of \emph{human understanding} and \emph{authority}.

Here's an example of the kind of proposal I want to see for ``natural'' mathematics. We could replace traditional authorship with \emph{co-ownership} of mathematical ideas. In this paradigm, any mathematician who demonstrates \emph{authoritative} understanding of a work --- the type you would expect of an author today --- would be entitled to claim co-ownership, even after publication. Some papers might have a few co-owners; others would have tens, or even hundreds. Journals would have the exciting, but challenging, role of establishing norms for validating understanding and maintaining the infrastructure of co-ownership. This would take time --- immense amounts of it. Explaining an entire paper in full detail to an appropriately skeptical audience often constitutes an entire graduate course. But if mathematicians aren’t writing papers anymore, we will have more time. That time should be returned to mathematics, in its most social and human forms. I think it sounds like fun.

\subsection*{For institutions}
Finally, there is work ahead for national and international mathematics institutions. In framing our work to the public, choosing the distribution of awards and grants, and organizing conferences and events, institutions have many tools to influence our path forward. They should use this power to make the case that mathematical knowledge happens between humans, not between computers. That means grants, awards, and conferences for promoting natural mathematics. We all know that mathematicians don't like politics, but a little political thinking would go a long way here. Institutions will act faster if their constituencies push them to action. You can be a part of moving them in the right direction. 

Natural mathematics can be an alternative to an artificial future. 
Some of these ideas will be good, and others will be bad; they’re new, and they need to be tested. Most importantly, we need more ideas for what natural mathematics should be. Let's build it together. Find your like-minded colleagues now, get your organizational meeting on the calendar, and let’s reclaim the future of our field.

\section*{Acknowledgments}
This essay was written in response to calls from Alper, Tao, and Venkatesh for visions for the future of mathematics \cite{MR5059809, taoicm, MR4726987}. While I disagree with some of their viewpoints, I am deeply indebted to their ideas, and grateful to them for beginning this public conversation. Juliet Glazer and Richard Schwartz provided valuable comments.

\bibliographystyle{abbrv}
\bibliography{bib}
\end{document}